\documentclass[10pt,a4paper]{amsart}

\usepackage{amsfonts,amsmath,amssymb,graphics}

\input xy
\xyoption{all}

\newtheorem{thm}{Theorem}[section]
\newtheorem{lemma}[thm]{Lemma}
\newtheorem{cor}[thm]{Corollary}
\newtheorem{pro}[thm]{Proposition}
\newtheorem{example}[thm]{Example}
\newtheorem{remark}[thm]{Remark}
\newtheorem{ddef}[thm]{Definition}

\def\tor{\mathop{\rm Tor}\nolimits}

\def\spec{\mathop{\rm Spec}}
\def\proj{\mathop{\rm Proj}}

\def\depth{\mathop{\rm depth}}
\def\reg{\mathop{\rm reg}\nolimits}
\def\Freg{\mathop{{\mathfrak f}{\rm reg}}\nolimits}

\def\a{{\alpha}}

\def\fin{\mathop{\rm end}}

\def\bideg{\mathop{\rm bideg}}
\def\cd{\mathop{\rm cd}\nolimits}
\def\ann{\mathop{\rm ann}}

\def\indeg{\mathop{\rm indeg}}
\def\bideg{\mathop{\rm bideg}}
\def\depth{\mathop{\rm depth}\nolimits}

\def\supp{\mathop{\rm Supp}}
\def\sym{\mathop{\rm Sym}\nolimits}

\def\E{{\mathcal E}}

\def\C{{\mathcal C}}

\def\MM{{\mathbb M}}

\def\OO{{\mathcal O}}

\def\M{{\mathrm M}}
\def\sF{{\mathcal F}}
\def\sI{{\mathcal I}}

\def\R{{\mathcal R}}

\def\a{{\alpha}}

\def\ip{{\mathfrak p}}
\def\ka{{\mathfrak K}}
\def\iq{{\mathfrak q}}

\def\Ip{{\mathfrak P}}

\def\ib{{\mathfrak b}}

\def\im{{\mathfrak m}}

\def\ra{{\rightarrow}}
\def\lra{{\longrightarrow}}
\def\fini{{$\quad\quad\Box$}}

\newcommand{\bd}{\begin{ddef}}
\newcommand{\ed}{\end{ddef}}
\newcommand{\bt}{\begin{thm}}
\newcommand{\et}{\end{thm}}
\newcommand{\bl}{\begin{lemma}}
\newcommand{\el}{\end{lemma}}
\newcommand{\bco}{\begin{cor}}
\newcommand{\eco}{\end{cor}}
\newcommand{\bp}{\begin{pro}}
\newcommand{\ep}{\end{pro}}
\newcommand{\bex}{\begin{example}}
\newcommand{\eex}{\end{example}}
\newcommand{\brm}{\begin{remark}}
\newcommand{\erm}{\end{remark}}
\newcommand{\bconj}{\begin{conj}}
\newcommand{\econj}{\end{conj}}

\newcommand{\fn}{\frak{n}}

\newcommand{\beqn}{\begin{eqnarray*}}
\newcommand{\eeqn}{\end{eqnarray*}}
\newcommand{\beq}{\begin{eqnarray}}
\newcommand{\eeq}{\end{eqnarray}}
\newcommand{\been}{\begin{enumerate}}
\newcommand{\eeen}{\end{enumerate}}

\begin{document}

\author[Chardin]{Marc Chardin}

\address{M. Chardin, Institut de Math\'ematiques de Jussieu,
CNRS \& UPMC,
4 place Jussieu,
F-75005, Paris, France.}
\email{chardin@math.jussieu.fr} 

\title[Powers of ideals]{Powers of ideals and the cohomology of  stalks and fibers of morphisms}

\maketitle

\begin{abstract}
We first provide here a very short  proof of a refinement of a theorem of Kodiyalam and Cutkosky, Herzog and Trung
on the regularity of powers of ideals. This result implies a conjecture of H\`a and generalizes a result of Eisenbud and Harris
concerning the case of ideals primary for the graded maximal ideal in a standard graded algebra over a field. It also implies a new result on the regularities of powers of ideal sheaves. We then compare 
the cohomology of the stalks and the cohomology of the fibers of a projective morphism to the effect of comparing the maximum over fibers and over stalks of the Castelnuovo-Mumford regularities of a
family of projective schemes.
\end{abstract}
\section{Introduction}

An important result of Kodiyalam and Cutkosky, Herzog and Trung states that the Castelnuovo-Mumford regularity
of the power $I^t$ of an ideal over a standard graded algebra is eventually a linear function in $t$. The leading term
of this function has been determined by Kodiyalam in his proof. 

This result was first obtained for standard graded algebras over a field, and later extended by Trung and Wang to 
standard graded algebras over a Noetherian ring. 

We first provide here a very short  proof of a refinement of this result

\bt
Let $A$ be a positively graded Noetherian algebra, $M\not= 0$ be a finitely generated graded $A$-module, 
$I$ be a graded $A$-ideal, and set 
$$
d:=\min\{ \mu\ \vert\ \exists p, \ (I_{\leq \mu})I^pM =I^{p+1}M\} .
$$
Then
$$
\lim_{t\ra \infty}(\fin (H^i_{A_+}(I^t M))+i-td)\in {\bf Z}\cup\{ -\infty\}
$$
exists for any $i$, and is at least equal to the initial degree of $M$ for some $i$.
\et

The end of a graded module $H$ is $\fin (H):=\sup\{ \mu \ \vert\ H_\mu \not= 0\}$ if $H\not= 0$ and $-\infty$ else.
Recall that for a graded $A$-module $N$, $\reg (N)=\max_i \{\fin (H^i_{A_+}(N))+i\}$. 

Very interesting examples showing a hectic behaviour of the value of $a^i(t):=\fin (H^i_{A_+}(I^t))$ as $t$ varies were given
by Cutkosky in \cite{Cu}. These examples point out that the existence of the limit quoted above do not imply that all
of the functions $a^i(t)$ are eventually linear functions of $t$. It only implies that at least one of them is eventually linear in $t$. For instance, in the examples given by Cutkosky, the limit in the Theorem is $-\infty$ for all $i\not= 0$.

More recently, Eisenbud and Harris proved that in the case of a standard graded algebra $A$ over a field, for a graded ideal
which is $A_+$-primary and generated in a single degree, the constant term in the linear function is the maximum
of the regularity of the fibers of the morphism defined by a set of minimal generators. In a recent preprint, Huy T\`ai  H\`a 
generalized this result by proving that, if an ideal is generated in a single degree $d$, a variant of the regularity (the $a^*$-invariant)
satisfies $a^* (I^t)=dt+a$ for $t\gg 0$, where $a$ can be expressed in terms of the maximum of the values of $a^*$ on the stalks of
the projection $\pi$ from the closure of the graph of the map defined by the generators to its image,  \cite[1.3]{Ha}. He conjectures that a 
similar result holds for the regularity.

In Theorem \ref{regpowgeo} we prove this conjecture of H\`a. More precisely, we show that the limit in the Theorem above is the
maximum of the end degree of the $i$-th local cohomology of the stalks of $\pi$, for ideals generated in a single degree. 
This holds for graded ideals in a Noetherian positively graded algebra.

An interesting, and perhaps surprising, consequence of this result is the following result on the limit of the regularity of saturation of powers, or equivalently of
powers of ideal sheaves, in a positively graded Noetherian algebra :

\bco
Let $I$ be a graded ideal generated in a single degree $d$. Then,
$$
\lim_{t\ra \infty}(\reg ((I^t)^{sat})-dt)
$$
exists and the following are equivalent :

(i) the limit is non negative,

(ii) the limit is not $-\infty$,

(iii) the projection $\pi$ from the closure of the graph of the  function defined by minimal generators of $I$ 
to its image admits a fiber of positive dimension.
\eco

This can be applied to ideals generated in degree at most $d$, replacing $I$ by $I_{\geq d}$. 

It gives a simple geometric criterion 
for an ideal $I$ generated in degree (at most) $d$ to satisfy $\reg ((I^t)^{sat})=dt+b$ for $t\gg 0$ : this holds if and only if there exists a subvariety $V$ of the closure of the graph that is contracted in its projection to  the closure of the image ({\it i.e.} $\dim (\pi (V))<\dim V$). A very simple example
is the following : in a polynomial ring in $n+1$ variables any graded ideal generated by $n$ forms of same degree $d$ satisfies
$\reg ((I^t)^{sat})=dt+b$ for $t\gg 0$, with $b\geq 0$. The same result holds if a reduction of the ideal is generated by at most $n$ 
elements (in other words if the analytic spread of $I$ is at most $n$).

The result of Eisenbud and Harris is stated in terms of regularity of fibers. For a finite morphism, there is no difference between 
the regularity of stalks and the regularity of fibers. This follows from the following result that is likely part of folklore, but that we didn't 
find in several of the classical references in the field :

\bl
Let $(R,\im ,k)$ be a Noetherian local ring,  $S:=R[X_1,\ldots  ,X_n]$ be a polynomial ring 
over $R$ with $\deg X_i>0$ and $\MM$ be a finitely generated graded $S$-module. 
Set $d:=\dim (\MM \otimes_R k )$. 
 Then  $H^{i}_{S_+}(\MM )=0$ for $i>d$ and  the natural graded map $H^{d}_{S_+}(\MM )\otimes_{R} k \lra H^{d}_{S_+}(\MM \otimes_{R} k) $ is an isomorphism.
\el

For morphisms that are not finite or flat, the situation is more subtle -- see Proposition \ref{CohFib}. We show that for
families of projective schemes that are close to being flat (if the Hilbert polynomial of any two fibers differ at most by a constant, in the standard graded situation),
the maximum of the regularities of stalks and the maximum of the regularities of fibers agree. Also the maximum regularity of stalks bounds
above the one for fibers under a weaker hypothesis. Putting this together provides a collection of results that covers the
results obtained in \cite{EH} and \cite{Ha}. See Theorem \ref{RegPow}.

To simplify the statements, we introduce the notion of regularity over a scheme, generalizing the usual notion of regularity
with reference to a polynomial extension of a ring. This is natural in our situation : the family of schemes given by the closure of the graph 
over the parameter space given by the closure of the image of our map, considered as a projective scheme,  is a key ingredient of this study. 
 
 This work was inspired by results of Huy T\`ai H\`a in  \cite{Ha} and of David Eisenbud and Joe Harris in \cite{EH}.  Bernd Ulrich made remarks on a very early version of some of these results and motivated my study of the difference between the 
 regularity of stalks and the regularity of fibers, and Joseph Oesterl\'e provided references concerning Lemma \ref{SuppTor}. It is
 my pleasure to thank  them for their contribution.
 
\section{Notations and general setup}

Let $R$ be a commutative ring and $S$ a polynomial ring over $R$ in finitely many variables.

If $S$ is ${\bf Z}$-graded, $R\subset S_0$, and $X_1,\ldots ,X_n$ are the variables with positive degrees, the \v Cech complex 
$\C^\bullet_{S_+} (M)$ with $\C^0_{S_+}(M)=M$ and $\C^i_{S_+}(M)=\oplus_{j_1<\cdots j_i}M_{X_{j_1}\cdots X_{j_i}}$ for $i>0$, is graded, whenever $M$ is a graded $S$-module.

There is an isomorphism $H^i_{S_+}(M)\simeq H^i (\C^\bullet_{S_+} (M))$ for all $i$, which is graded if $M$ is so.
One then defines two invariants attached to such a graded $S$-module $M$ :
$$
a^i (M):=\sup \{ \mu\ \vert\ H^i_{S_+}(M)_\mu \not= 0\}
$$
if $H^i_{S_+}(M)\not= 0$ and $a^i (M):=-\infty$ else, and
$$
b_j (M):=\sup \{ \mu\ \vert\ \tor_j^{S}(M,S/S_+)_\mu \not= 0\}
$$
if $\tor_j^{S}(M,S/S_+)\not= 0$ and $b_j (M):=-\infty$ else. Notice that $a^i (M)=-\infty$ for $i>n$ and 
$b_j (M)=-\infty$ for $j>n$. The Castelnuovo-Mumford regularity of a graded $S$-module
$M$ is then defined as 
$$
\reg (M):=\max_{i}\{ a^i (M)+i\} =\max_j \{b_j (M)-j\} +n-\sigma
$$
where $\sigma $ is the sum of the degrees of the variables with positive degrees. Other options are possible, 
in particular when $S$ is not standard graded (when $\sigma \not= n$).  
Another related invariant is 
$$
a^* (M):=\max_{i}\{ a^i (M)\} =\max_j \{b_j (M)\} -\sigma .
$$

The following  classical result is usually stated for positive grading.

\bt\label{Zgraded} Let $S$ be a Noetherian ${\bf Z}$-graded algebra and $M$ be a finitely
generated graded $S$-module. Then, for any $i$, 

(i) $a^i (M) \in \{ -\infty\} \cup {\bf Z}$,

(ii) the $S_0$-module $H^i_{S_+}(M)_\mu$ is finitely generated for any $\mu \in {\bf Z}$.
\et

{\it Proof.} As $S$ is Noetherian, $S$ is an epimorphic image of a polynomial ring $S'$ over $S_0$ 
by a graded morphism. Considering $M$ as $S'$-module one has $H^i_{S_+}(M)\simeq  H^i_{S'_+}(M)$ via
the natural induced map, so that we may replace $S$ by $S'$ and assume that 
$$
S=S_0 [Y_1,\ldots ,Y_m,X_1,\ldots ,X_n]
$$ 
with $\deg Y_i\leq -1$ and $\deg X_j\geq 1$ for all $i$ and $j$. We recall that 
$H^i_{S_+}(S)=0$ for $i<n$ and $H^n_{S_+}(S)=(X_1\cdots X_n)^{-1}S_0 [Y_1,\ldots ,Y_m,X_1^{-1},\ldots ,X_n^{-1}]$, and notice that $H^n_{S_+}(S)_\mu$ is a finitely generated free $S_0$-module for any 
$\mu$.

Let $F_\bullet$ be a graded free $S$-resolution of $M$ with $F_i$ finitely generated. Both spectral sequences
associated to the double complex $\C_{S_+}^\bullet F_\bullet$ degenerate at step 2 and provide graded 
isomorphisms :
$$
H^i_{S_+}(M)\simeq H_{n-i}(H^n_{S_+}(F_\bullet )),
$$
which shows that $H^i_{S_+}(M)_\mu$ is a subquotient of $H^n_{S_+}(F_{n-i})_\mu$,
hence a finitely generated $S_0$-module which is zero in degrees $>-n+b_{n-i}$, where 
$b_j$ is the highest degree of a basis element of $F_j$ over $S$.\fini

\section{Regularity over a scheme.}



Local cohomology and the torsion functor commute with localization on the base $R$, providing
natural graded isomorphisms for a graded $S$-module $M$ :
$$
H^i_{(S\otimes_R R_\ip)_+}(M\otimes_R R_\ip)\simeq H^i_{S_+}(M)\otimes_R R_\ip
$$
and
$$
\tor_i^{S\otimes_R R_\ip}(M\otimes_R R_\ip , R_\ip)\simeq \tor_i^{S}(M,R)\otimes_R R_\ip .
$$
Hence
$a^i (M)=\sup_{\ip \in \spec (R)}a^i (M\otimes_R R_\ip )$ and  $b_j (M)=\sup_{\ip \in \spec (R)}b_j (M\otimes_R R_\ip )$. It follows that the regularity is a local notion on $R$ :
$$
\reg (M)=\sup_{\ip \in \spec (R)}\reg (M\otimes_R R_\ip ).
$$
These supremums are maximums whenever $\reg (M)<+\infty$, for instance if $R$ is Noetherian and $M$ is finitely generated. The same holds for $a^* (M)$.

Extending this definition to the case where the base is a scheme is  natural and is given in the following definition.

\bd\label{RegSch} Let $Y$ be a scheme, $\E$ be a locally free $\OO_Y$-module of finite rank, $\sF$ be a 
graded sheaf  of $\sym_Y (\E )$-modules. Then
$$
a^i (\sF ):=\sup_{y\in Y}a^i (\sF \otimes_{\OO_Y}  \OO_{Y,y}), \quad \reg (\sF ):=\max_{i}\{ a^i (\sF )+i\} .
$$

If $\E$ is free, $\sym_Y (\E )={\OO_Y}[X_1,\ldots ,X_n]$, and the definition of regularity above makes sense
for non standard grading. 

A closed subscheme $Z$ of  $\proj (\sym_Y (\E ))$ correspond to a unique
graded $\sym_Y (\E )$-ideal sheaf $\sI_Z$ saturated with respect to $\sym_Y (\E )_+$. We set 
$$
a^i (Z):=\sup_{y\in Y}a^i ({\OO_{Y,y}}[X_0,\ldots ,X_n]/(\sI_Z \otimes_{\OO_Y} \OO_{Y,y}))
$$
(notice that $a^0 (Z)=-\infty$) and $\reg (Z):=\max_{i}\{ a^i (Z)+i\}$.
\ed


The following Proposition is immediate from the definition and the corresponding results 
over an affine scheme. 

\bp
Assume $Y$ is Noetherian, $\E$ is a locally free coherent sheaf on $Y$ and $\sF \not= 0$ is a coherent
graded sheaf  of $\sym_Y (\E )$-modules. Then $\reg (\sF )\in {\bf Z}$.
If $Z\not= \emptyset$ is a closed subscheme of ${\bf P}^{n-1}_Y$,
then $\reg (Z)\geq 0$.
\ep

\section{First result on cohomology of powers.}

We now prove the first statement of our text on cohomology of powers of ideals. It refines earlier results on the regularity of powers (\cite{Ko}, \cite{CHT} and \cite{TW}). The 
argument is based on Theorem \ref{Zgraded} applied to a Rees algebra  and a lemma due to Kodiyalam.

\bt\label{AB1}
Let $A$ be a positively graded Noetherian algebra, $M\not= 0$ be a finitely generated graded $A$-module, 
$I$ be a graded $A$-ideal, and set 
$$
d:=\min\{ \mu\ \vert\ \exists p, \ (I_{\leq \mu})I^pM =I^{p+1}M\} .
$$
Then
$$
\lim_{t\ra \infty}(a^i (I^t M)+i-td)\in {\bf Z}\cup \{ -\infty\}
$$
exists for any $i$, and is at least equal to $\indeg (M)$ for some $i$.
\et

{\it Proof.} Set $J:=I_{\leq d}$ and write $J=(g_1,\ldots ,g_s)$ with $\deg g_i =d$ for $1\leq i\leq m$ and
$\deg g_i<d$ else. Let
$$
\R_J :=\oplus_{t\geq 0} J(d)^t=\oplus_{t\geq 0} J^t(td)\quad \hbox{and}\quad 
\R_I :=\oplus_{t\geq 0} I(d)^t=\oplus_{t\geq 0} I^t(td),
$$
and $S_0:=A_0 [T_1,\ldots ,T_m]$, $S:=S_0 [T_{m+1},\ldots ,T_s,X_1,\ldots ,X_n]$, with $\deg (T_i):=
\deg (g_i)-d$. Setting $\bideg (T_i):=(\deg (T_i) ,1)$ and $\bideg (X_j):=(\deg (X_j) ,0)$, one has 
$J(d)=(\R_J)_{0,1}$, hence a  bigraded onto map
$$
\xymatrix@R=0pt{
S\ar[r]&\R_J\\
T_i\ar@{|->}[r]&g_i \in J(d).\\
}
$$
As $M\R_I$ is finite over $\R_J$ 
according to the definition of $d$, the bigraded embedding $\R_J \lra \R_I$ makes $M\R_I$ a finitely generated bigraded $S$-module. 

The equality of graded $A$-modules $H^i_{S_+}(M\R_I )_{(* ,t )}=H^i_{A_+}(M\R_I )_{(* ,t )}$ 
shows that
$$
H^i_{S_+}(M\R_I )_{(\mu ,t )}=H^i_{A_+}((M\R_I )_{(*,t)})_{\mu}=H^i_{A_+}(MI^t)_{\mu +td}.
$$

By Theorem \ref{Zgraded} (i), $a^i (M\R_I )<+\infty$ and the above equalities show that $a^i (MI^t)\leq td+a^i (M\R_I )$, 
and that equality holds for some $t$. 

Furthermore, Theorem \ref{Zgraded} (ii) shows that $K_{i,\mu}:=H^i_{S_+}(M\R_I )_{(\mu ,* )}$ is a finitely 
generated graded $S_0$-module (for the standard grading $\deg (T_i)=1$). It follows that $H^i_{S_+}(M\R_I )_{(\mu ,t )}=0$ for $t\gg 0$ if and
only if $K_{i,\mu}$ is annihilated by a power of $\fn :=(T_1,\ldots ,T_m)$. Hence
$$
\lim_{t\ra +\infty} (a^i (MI^t)-td)=-\infty
$$
if $K_{i,\mu}$ is annihilated by a power of $\fn$ for every $\mu \leq a^i (M\R_I )$,
and else
$$
\lim_{t\ra +\infty} (a^i (MI^t)-td)=\max \{ \mu\ \vert\ K_{i,\mu}\not= H^0_\fn (K_{i,\mu}) \} .
$$

As $\reg (MI^t )\geq \fin (MI^t/R_+MI^t)$, the last claim follows from the next lemma, due to Kodiyalam.
\fini

\bl
With the hypotheses of Theorem \ref{AB1}, 
$$
\fin (MI^t/A_+MI^t)\geq \indeg (M)+ td, \quad \forall t.
$$
\el

{\it Proof.} The proof goes along the same lines as in [Ko, the proof of Proposition 4]. Notice that
the needed graded version of Nakayama's lemma applies.\fini

\section{Cohomology of powers and cohomology of stalks}

The following result is a more elaborated, and more technical,  version of Theorem \ref{AB1} that essentially follows from 
its proof. It implies a conjecture of H\`a on the regularity of powers of ideals, 
and refines the main result in \cite{Ha}. We will see later that, combined with a result on the regularity of
stalks and fibers of a morphism,  it also implies the result of Eisenbud and 
Harris in \cite{EH}. 

\bp\label{AB2}
Let $A$ be a positively graded Noetherian algebra, $M$ be a finitely generated graded $A$-module, 
$I$ be a graded $A$-ideal and $J\subseteq I$ be a graded ideal such that  $JI^pM=I^{p+1}M$ for some $p$. 

Assume that the ideal $J$ is generated by $r$ forms $f_1,\ldots ,f_r$ of respective degrees $d_1=\cdots =d_m>d_{m+1}\geq \cdots \geq d_r$. Set $d:=d_1$, $\deg (T_i):=\deg (f_i)-d$, $\bideg (T_i):=(\deg (T_i),1)$ and 
$\bideg (a):=(\deg (a),0)$ for $a\in A$. Consider the natural bigraded morphism of  bigraded 
$A_0$-algebras 
$$
\xymatrix{
S:=A[T_1,\ldots ,T_r]\ar^(.4){\psi}[r]&\R_I:=\oplus_{t\geq 0} I(d)^t=\oplus_{t\geq 0} I^t(dt),\\
}
$$
sending $T_i$ to $f_i$ and the bigraded map of $S$-modules :
$$
\xymatrix{
M[T_1,\ldots ,T_r]\ar^(.45){1_M \otimes_A \psi }[rr]&&M\R_I:=\oplus_{t\geq 0}M I^t(dt).\\
}
$$
Let  $B:=A_0[T_1,\ldots ,T_m]$ and $B':=B/\ann_B (\ker (1_M \otimes_A \psi ))$. 

Then,
$$
\lim_{t\ra +\infty} (a^i (MI^t)-td)=\max_{\iq \in \proj (B')}\{ a^i (M\R_I \otimes_{B'}B'_\iq )\}.
$$
\ep

{\it Proof.} 
First remark that in the proof of Theorem \ref{AB1} we only need the equality $JI^pM=I^{p+1}M$ for some $p$ (as a consequence, for all $p$
big enough).  We have shown there that 
$$
(*)\quad \lim_{t\ra +\infty} (a^i (MI^t)-td)=-\infty ,
$$
if and only if the finitely generated $B$-module $H^i_{S_+}(M\R_I )_{(\mu ,* )}$ is supported in  $V(T_1,\ldots ,T_m)$ for any $\mu$. As local cohomology commutes with flat base change
and elements in $B$ have degree 0,
$$
H^i_{S_+}(M\R_I )_{(\mu ,* )}\otimes_{B'}B'_\iq =H^i_{S_+}(M\R_I \otimes_{B'}B'_\iq )_{(\mu ,* )}
$$
hence  $(*)$ holds if and only if $H^i_{S_+}(M\R_I \otimes_{B'}B'_\iq )=0$ for any $\iq \in \proj (B')$. On the other hand, if this does not hold, there exists 
$\mu_0$ maximum  such that $H^i_{S_+}(M\R_I )_{(\mu_0 ,* )}$ is not supported in  $V(T_1,\ldots ,T_m)$,
and choosing $\iq\in \proj (B')\cap \supp (H^i_{S_+}(M\R_I )_{(\mu_0 ,* )})$ shows that both members in the asserted equality are equal to $\mu_0$.
\fini

\brm
In the above Proposition, as well as in other places in this text, we localize at homogeneous 
primes $\iq \in \proj (C)$ for some standard graded algebra $C$, in other words at graded prime
ideals that do not contain $C_+$. We may as well replace these localization by the degree zero part
of the localization at such a prime ideal, usually denoted by $C_{(\iq )}$ : the multiplication by
an element $\ell \in C_1\setminus \iq$ induces an isomorphism  $(C_{\iq })_\mu \simeq (C_{\iq })_{\mu +1}$ for any $\mu$. 
Hence, for any $C$-module $M$, $M\otimes_C C_\iq =0$ if and only if $M\otimes_C C_{(\iq )} =0$.
\erm

In the equal degree case,  the following corollary, that we state in a more geometric fashion, implies the conjecture of H\`a in \cite{Ha}.

\bt\label{regpowgeo}
Let $A:=A_0[x_0,\ldots ,x_n]$ be a positively graded Noetherian algebra and $I$ be a graded $A$-ideal generated
by $m+1$ forms of degree $d$. Set $Y:=\spec (A_0)$ and $X:=\proj (A/I)\subset \proj (A)\subseteq  \tilde{\bf P}_Y^n$. 
Let $\phi :\tilde{\bf P}_Y^n\setminus X\lra {\bf P}^{m}_Y$ be the
corresponding rational map, $W$ be the closure of the image of $\phi$, and
$$
\Gamma \subset \tilde{\bf P}^n_{W}\subseteq \tilde{\bf P}^n_{{\bf P}^{m}_Y}=\tilde{\bf P}_Y^n \times_Y {\bf P}^{m}_Y
$$ be the closure of the graph of $\phi$. Let  
$
\pi : \Gamma \lra W
$ 
be the projection induced by the natural map $ \tilde{\bf P}^n_{{\bf P}^{m}_Y}\lra {\bf P}^{m}_Y$. Then
 $$
 \lim_{t\ra +\infty}(a^i (I^t)-dt)=a^i (\Gamma ).
 $$
\et

{\it Proof.} Choose $J:=I$ and $M:=A$ in Proposition \ref{AB2}.
The equality $\lim_{t\ra +\infty}(a^i (I^t)-dt)=a^i (\Gamma )$ directly follows from the conclusion of Proposition \ref{AB2} according the definition of 
$a^i(\Gamma)$ for $\Gamma \subset \tilde{\bf P}^n_{W}$ given in Definition \ref{RegSch}.\fini

\section{Cohomology of stalks and cohomology of fibers}

We will now compare the cohomology of stalks and of fibers of a projective morphism, 
in order to compare their Castelnuovo-Mumford regularities. It will need results on the support
of Tor modules.  These are likely part of folklore. However, we included a proof as we did not find a reference
that properly fits our exact need.

\bl\label{SuppTor}
Let $R \ra S$ be a  homomorphism of Noetherian rings, $\MM $ be a finitely generated $S$-module and $N$ be a finitely generated $R$-module. 

Then the
$S$-modules $\tor_q^R (\MM , N )$ are finitely generated over $S$ and

(i) $\supp_S (\tor_q^R (\MM , N ))\subseteq \supp_S (\MM \otimes_R N)$ for any $q$, 

(ii) if further $(R ,\im)$ is local,  $S=R[X_1,\ldots, X_n]$, with $\deg X_i>0$ and $\MM $
is a graded $S$-module,
then $\supp_S (\tor_q^R (\MM , R/\im ))\subseteq \supp_S (\tor_1^R (\MM , R/\im ))$ 
for any $q\geq 1$.
\el

{\it Proof.} First the modules $\tor_q^R (\MM , N )$ are finitely generated over $S$ by \cite[X \S 6 N$^{\circ}$4 Corollaire]{BA}. Second, $\supp_S (\MM \otimes_R N)=\supp_S (\MM )\cap \varphi^{-1}(\supp_R (N))$, where
$\varphi : \spec (S)\lra \spec (R)$ is the natural map induced by $R \lra S$, by \cite[II \S 4 N$^{\circ}$4, Proposition 18 \& Proposition 19]{BAC}, since 
$\MM \otimes_R N =\MM \otimes_S (N\otimes_R S)$.
For $\Ip\in \spec (S)$, set $\ip :=\varphi (\Ip )$. Then   $\tor_q^R (\MM , N )_{\Ip}=\tor_q^{R_{\ip}} (\MM _\Ip , 
N_\ip )$ vanishes if either  $\MM _\Ip =0$ or $N_\ip =0$. 

For (ii), we can reduce to the case of a local morphism by localizing $S$ and $\MM $ at $\im +S_+$. 
In this local situation, $\tor_1^R (\MM , R/\im )=0$ if and only if $\MM $ is $A$-flat by \cite[Lemme 58]{An},
which proves our claim by localization at primes $\Ip$ such that $\varphi (\Ip )=\im$.\fini

Let $R$ be a commutative ring, $N$ be a $R$-module, $S:=R[X_1,\ldots  ,X_n]$ a polynomial ring
over $R$ and $\MM $ be a graded $S$-module. For a $S$-module $\MM$, we will 
denote by $\cd_{S_+} (\MM )$ the cohomological dimension of $\MM$ with respect to $S_+$, which is the maximal index $i$ such that $H^i_{S_+}(\MM )\not= 0$ (and $-\infty$ if all these local cohomology groups are 0).  
The following lemma is a natural way for 
comparing cohomology of stalks to cohomology of fibers.

\bl\label{CompCoh}
There are two converging spectral sequences of graded $S$-modules with same abutment $H^\bullet$ and with respective second terms
$$
{_2 'E}^p_q=H^p_{S_+} (\tor_q^R(\MM ,N ))\Rightarrow H^{p-q}
$$
and
$$
{_2 ''E}^p_q=\tor_q^R(H^p_{S_+}(\MM ),N)\Rightarrow H^{p-q}.
$$
Let $d :=\max \{ i\ \vert \ H^i_{S_+}(\MM \otimes_R N)\not= 0\}$.  If $R$ is Noetherian, $N$ is finitely generated over 
$R$ and $\MM $ is finitely generated over $S$, then
$$
H^{d}_{S_+}(\MM \otimes_R N)\simeq H^{d}_{S_+}(\MM )\otimes_R N 
$$
and $\tor_q^R(H^{i}_{S_+}(\MM ), N) =H^i_{S_+}(\tor_q^R (\MM , N ))=0$ for any $q$ if $i>d$.
\el
{\it Proof.} Let $F_\bullet$ be a free $R$-resolution of $N$. Consider the 
double complex $\C^\bullet_{S_+}(\MM \otimes_ R F_\bullet )=\C^\bullet_{S_+}(\MM )\otimes_ R F_\bullet $, totalizing to 
$T^\bullet$ with $T^i=\oplus_{p-q=i}\C^p_{S_+}(\MM )\otimes_ R F_q $. It gives rise to two spectral 
sequences abuting to the homology $H^\bullet$ of $T^\bullet$. 

One has  first terms $\C^p_{S_+}(\tor_q^R (\MM , N ))$ and second terms $H^p_{S_+}(\tor_q^R (\MM , N ))$.

The other spectral sequence has first terms  $H^p_{S_+}(\MM )\otimes_ R F_q$
and second terms $\tor_q^R(H^p_{S_+}(\MM ),N)$. 
It provides the quoted spectral sequences.

Recall that if $P$ is a finitely presented $S$-module, one has 
$\cd_{S_+}(P')\leq \cd_{S_+}(P)$ whenever
$\supp (P')\subseteq \supp (P)$. This is proved in \cite[2.2]{DNT} under the assumption that
$S$ is Noetherian and $P'$ is finitely generated, which is enough for our purpose. 

By Lemma \ref{SuppTor} (i), $\supp (\tor_q^R (\MM , N ))\subseteq \supp (\MM \otimes_R N)$ for any $q$, 
which implies that
$H^i_{S_+}(\tor_q^R (\MM , N ))=0$ for any $q$ if $i>d$. It follows that 
$H^{d}=H^{d}_{S_+}(\MM \otimes_R N)$ and $H^{i}=0$ for $i>d$. 

On the other hand, choose $i$ maximal such that $H^i_{S_+}(\MM )\otimes_R N\not= 0$. 
Then $\tor_q^R(H^p_{S_+}(\MM ),N)=0$ for any $q$ if $p>i$, because $H^p_{S_+}(\MM )_\mu$
is a finitely generated $R$-module for every $\mu$,
and hence 
$H^{i}=H^{i}_{S_+}(\MM )\otimes_R N\not= 0$ and $H^{j}=0$ for $j>i$. The conclusion follows.
\fini

The following statement extends a classical result on the cohomolgy of fibers in
a flat family (see for instance  \cite[III 9.3]{Har}). The hypothesis on the cohomological dimension of Tor
modules that appears in (ii) will be connected to the variation of the Hilbert polynomial of fibers in the
corresponding family of sheaves in Proposition \ref{TorHP} ; it is a weakening of the flatness condition
for this family.

\bp\label{CohFib}
Let $(R,\im ,k)$ be a Noetherian local ring,  $S:=R[X_1,\ldots  ,X_n]$ be a polynomial ring 
over $R$, with $\deg X_i>0$ for all $i$, and $\MM$ be a finitely generated graded $S$-module. Set $\M:=\MM \otimes_R k$ and $d:=\dim \M$. 
 Then one has,\\

(i) The natural graded map $H^{d}_{S_+}(\MM )\otimes_{R} k \lra H^{d}_{S_+}(\M )$ is an isomorphism and $d =\max\{ i\ \vert \ H^i_{S_+}(\MM ) \not= 0\} $.
In particular, 
$$
a^{{d}}(\MM )=a^{{d}}(\M ) \in {\bf Z}.
$$\\
(ii) For any integers $\mu$ and $\ell$,  if $\cd_{S_+}(\tor_1^{R}(\MM ,k ))\leq \ell +1$ then
$$
\{ H^i_{S_+}(\MM )_\mu =0,\forall i\geq \ell\} \Rightarrow \{ H^i_{S_+}(\M )_\mu =0,\forall i\geq \ell\} ,
$$
and both conditions are equivalent if $\cd_{S_+}(\tor_1^{R}(\MM ,k ))\leq \ell$. 
In particular,  $\reg (\M ) \leq \reg (\MM  )$  if $\cd_{S_+}(\tor_1^{R}(\MM ,k ))\leq 1$ 
and equality holds if $\depth_{S_+}(\MM )>0$.\\
\ep

{\it Proof.} We  consider the two spectral sequences in Lemma \ref{CompCoh}, 

$$
{_2 'E}^p_q=H^p_{S_+} (\tor_q^{R } (\MM ,k ))\Rightarrow H^{p-q}
$$
and
$$
{_2 ''E}^p_q=\tor_q^{R }(H^p_{S_+}(\MM ),k )\Rightarrow H^{p-q}.
$$

Let $B:=k [X_1,\ldots ,X_n]$. 
The module $\tor_q^{R } (\MM ,k )$ is a $R [X_1,\ldots ,X_n]$-module of finite type,
annihilated by $\im$ and $\ann_{S} (\MM )$. Hence $\M$ is a graded $B$-module of finite type and $\tor_q^{R } (\MM ,k )$ is a graded 
$(B/\ann_{B} (\M))$-module of  finite type, for any $q$.

Notice that $d=\cd_{S_+}(\M)=\cd_{B_+}(\M)$. It follows that ${_2 'E}^p_q=0$ if $p>d $, and ${_2 'E}^{d}_0\not= 0$.

By Lemma \ref{CompCoh},
${_2 ''E}^p_q=0$ for all $q$ if $p>d $, in particular
$H^p_{S_+}(\MM )_\mu \otimes_{R} k =0$ for any $\mu$ if $p>d $. Hence 
$H^p_{S_+}(\MM )_\mu =0$ for any $\mu$ if $p>d $. In other words,
$H^p_{S_+}(\MM )=0$ for any $p>d $.

The same lemma shows that $H^{d}_{S_+}(\M )=H^{d}_{S_+}(\MM )\otimes_{R} k$, and finishes the proof of (i).

For (ii), let $\mu$ be an integer. We prove the result by descending 
induction on $\ell$ from the case $\ell =d$, which we already proved. 

Assume 
the results hold for $\ell +1$. Recall that, for any $p$, the map ${_r 'd}^{p-r}_{1-r}: {_r 'E}^{p-r}_{1-r} \lra {_r 'E}^{p}_{0}$ is zero for  $r\geq 2$ and
the map ${_r ''d}^{p}_{0}: {_r ''E}^{p}_{0} \lra {_r ''E}^{p+1-r}_{-r}$ is the zero map for $r\geq 1$.

If $H^i_{S_+}(\MM )_\mu =0,\forall i\geq \ell$, then 
$({_2 ''E}^p_q)_\mu =0$ for $p\geq \ell$ and all $q$. As ${_2 ''E}^p_q =0$ for $q<0$, it follows that $({_2 ''E}^p_q)_\mu =0$
if $p-q\geq \ell$. 

If $\cd_{S_+}(\tor_1^{R}(\MM ,k ))\leq \ell +1$ then ${_2 'E}^p_q=0$ for $p\geq \ell +2$ and 
$q>0$ by Lemma \ref{SuppTor} (ii), in particular the map
$$
({_r 'd}^{\ell }_{0})_\mu: ({_r 'E}^{\ell}_{0})_\mu \lra ({_r 'E}^{\ell +r}_{r-1})_\mu
$$
is the zero map  for any $r\geq 2$, hence 
$H^\ell_{S_+}(\M )_\mu =({_2 'E}^\ell_0)_\mu =({_\infty 'E}^\ell_0)_\mu =0$ as claimed.

For the reverse implication, the hypothesis implies that ${_2 'E}^p_q=0$ if $q\geq 1$ and $p\geq \ell +1$ by
Lemma \ref{SuppTor} (ii).
Hence $({_2 'E}^p_q)_\mu=0$ for $p-q\geq \ell$ if $H^\ell_{S_+}(\M )_\mu =0$. By induction hypothesis, $H^p_{S_+}(\MM )_\mu \otimes_{R} k
=0$ for $p\geq \ell +1$. Hence $({_2 'E}^p_q)_\mu =\tor_q^{R }(H^p_{S_+}(\MM )_\mu ,k )=0$ for $p\geq \ell +1$ and all $q$. It implies that $H^\ell_{S_+}(\MM )_\mu \otimes_{R} k =({_\infty ''E}^\ell_0)_\mu =0$, 
and proves the claimed equivalence.

Finally, recall that $H^i_{S_+}(\MM )=0$ for $i<\depth_{S_+}(\MM )$.\fini\medskip

\brm
Notice that  $\reg (\M ) \leq \reg (\MM  )$ does not hold without the hypothesis $\cd_{S_+}(\tor_1^{R}(\MM ,k ))\leq 1$.
To see this, consider generic polynomials of some given degrees $d_1,\ldots ,d_r$ : $P_i:=\sum_{\vert \a \vert =d_i}
U_{i,\a} X^\a \in k[U_{i,\a}] [X_1,\ldots ,X_n]$, with $r\leq n$ and a specialization map $\phi : k[U_{i,\a}]\lra  k$ to the field $k$ with kernel $\im$.
Set $R:=k[U_{i,\a}]_\im$. As the $P_i$'s form a regular sequence in $k[U_{i,\a}] [X_1,\ldots ,X_n]$,
they also form one in $S:=R[X_1,\ldots ,X_n]$ and show that $\MM :=S/(P_1,\ldots ,P_r)$ has regularity
$d_1+\cdots +d_r-r$. On the other hand,  the regularity of $M=k[X_1,\ldots ,X_n]/(\phi (P_1),\ldots ,\phi (P_r))$,
need not be bounded by $d_1+\cdots +d_r-r$. 

For instance, with $n=4$ and $r=3$,  take $\phi (P_1):=X_1^{d-1}X_2-X_3^{d-1}X_4$, $\phi (P_2):=X_2^d$
and $\phi (P_3):=X_4^d$ (over any field). Then one has $\reg (M)=d^2-2$ for $d\geq 3$ (see \cite[1.13.6]{Ch}) which is bigger than $\reg (\MM)=3d-3$, and $\cd_{S_+}(\tor_1^{R}(\MM ,k ))=2$.
\erm

\brm
In the other direction, it may of course be that $\reg (\MM )>\reg (M )$. If for instance $(R,\pi ,k)$ is
a DVR, one may take $\MM :=R[X]/(\pi X^d)$, so that $\reg (\MM )=d-1$ and $\reg (M)=0$, with 
$\cd_{S_+}(\tor_1^{R}(\MM ,k ))=1$.

More interesting is the example $R:={\bf Q}[a,b]$, $\im :=(a,b)$ and 
$$
\MM :=\sym_R(\im^3)=R[X_1,\ldots ,X_4]/(bX_1-aX_2,bX_2-aX_3,bX_3-aX_4).
$$
Then for any morphism from $R$ to a field $k$, $\reg (\MM\otimes_R k)=0$,
while $\reg (\MM )=1$. 

Similar examples arises from the symmetric algebra of  other ideals that are not generated by a proper sequence.
\erm

The characterization of flatness in terms of the constancy of the Hilbert polynomial of fibers extends as follows.


\bl\label{TorHP}
Let $p$ be an integer. In the setting of Proposition \ref{CohFib}, assume that $R$ is reduced and $S$ is standard graded. Then the following are equivalent :

(i) $\dim (\tor_1^{R}(\MM ,k ))\leq p$,

(ii) The Hilbert polynomials of $\MM \otimes_{R} k$ and $\MM \otimes_{R} (R_\ip /\ip R_\ip)$ differ at most by 
a polynomial of degree $<p$, for any $\ip\in \spec (R)$.
\el

{\it Proof.} 
We induct on $p$. The result is  standard  when $p=0$, see for instance \cite[III 9.9]{Har} or \cite[Ex. 20.14]{Ei}. 

Assume (i) and (ii) are equivalent for $p-1\geq 0$, for any Noetherian local domain, standard graded polynomial ring over it and graded module of finite type. 

Set $K:=R_\ip /\ip R_\ip$, $M_K:=\MM\otimes_R K$, $B:=k[X_1,\ldots ,X_n]$ and $C:=K[X_1,\ldots ,X_n]$. Consider variables $U_1,\ldots ,U_n$ (of degree 0) and let $\ell :=U_1X_1+\cdots +U_nX_n$. By the Dedekind-Mertens lemma,
$$
\begin{array}{l}
\hbox{(a)}\quad\quad\ker (\!\! \xymatrix{\MM [U]\ar^(.4){\times \ell}[r]&\MM [U](1)\\}\!\! )\subseteq H^0_{S_+}(\MM )[U],\\
\hbox{(b)}\quad\quad\ker (\!\! \xymatrix{M [U]\ar^(.4){\times \ell}[r]&M[U](1)\\}\!\! )\subseteq H^0_{B_+}(M)[U],\\
\hbox{(c)}\quad\quad\ker (\!\! \xymatrix{M_K [U]\ar^(.4){\times \ell}[r]&M_K [U](1)\\}\!\! )\subseteq H^0_{C_+}(M_K)[U],\\
\hbox{(d)}\quad\quad\ker (\!\! \xymatrix{\tor_1^{R}(\MM ,k ) [U]\ar^(.45){\times \ell}[r]&\tor_1^{R}(\MM ,k )[U](1)\\}\!\! )\subseteq H^0_{B_+}
(\tor_1^{R}(\MM ,k ))[U].\\
\end{array}
$$
Let $R':=R(U)$ be obtained from $R[U]$ by inverting all polynomials whose coefficient ideal is the unit ideal, and denote by $N'$ the extension of
scalars from $R$ to $R'$ for the module $N$. Recall that $R(U)$ is local with maximal ideal $\im R(U)$, residue field $k'=k(U)$ and that $K'=K(U)$ --see 
for instance [Na, p. 17]. As 
the zero local cohomology modules above vanish in high degrees, (b) and (c)
show that $\MM' /\ell \MM'$ satisfies condition (ii) of the Lemma for $p-1$, $R'$ and $R'[X_1,\ldots ,X_n]$. Now (a) and (d) provide an exact sequence for $\mu \gg 0$:
$$
\xymatrix@C=18pt{0\ar[r]&\tor_1^{R}(\MM ' ,k')_{\mu -1} \ar^(.5){\times \ell}[r]&\tor_1^{R'}(\MM ',k' )_{\mu } \ar[r]&\tor_1^{R'}(\MM '/\ell \MM ',k' ) _{\mu }\ar[r]&0\\}
$$
which shows in particular that 
$$
\dim \tor_1^{R'}(\MM '/\ell \MM ',k' )=\dim \tor_1^{R'}(\MM ',k' )-1=\dim \tor_1^{R}(\MM ,k)-1,
$$ 
if $\dim \tor_1^{R}(\MM ,k)$ is positive, and proves our claim by induction.
\fini\medskip

\brm

If the grading is not standard, a quasi-polynomial is attached to any finitely generated graded module, and in Lemma \ref{TorHP}
property  (ii) should be replaced by the following :

(ii) the difference between the quasi-polynomials of $\MM \otimes_{R} k$ and $\MM \otimes_{R} (R_\ip /\ip R_\ip )$ is a 
quasi-polynomial of degree $<p$ for any $\ip \in \spec (R)$.

The degree of a quasi-polynomial is the highest degree of the polynomials that defines it. The proof of \cite[III 9.9]{Har} extends to
this case when $p=0$, and our proof extends after a slight modification : in the proof of (ii)$\Rightarrow$(i), one should take
 $\ell :=U_1X_1^{w/w_1}+\cdots +
U_nX_n^{w/w_n}$, where $w_i:=\deg (X_i)$ and $w:=\hbox{lcm} (w_1,\ldots ,w_n)$. 
\erm

The local statement of Lemma \ref{TorHP} implies a global statement, by comparing Hilbert functions at generic points of
the components and at closed points. We state it below in a ring theoretic form. 

\bp\label{FlatnessCrit}
Let $p$ be an integer, $R$ be a reduced commutative ring, $S$ be a Noetherian positively graded polynomial ring over $R$ and
$\MM$ be a finitely generated graded $S$-module. Then the following are equivalent :

(i) $H^i_{S_+}(\tor_1^R(\MM ,R/\im ))=0$, for all $i>p$ and $\im$ maximal in $\spec (R)$,

(ii) for any two ideals $\ip\subset \iq$ in $\spec (R)$, the quasi-polynomials of $\MM \otimes_R R/\ip$ and $\MM \otimes_R R/\iq$
differ by a quasi-polynomial of degree $<p$, 

(iii) over a connected component of $\spec (R)$ the quasi-polynomials of two fibers differ 
by a quasi-polynomial of degree $<p$.
\ep

In parallel to the definition of the regularity over a scheme, we define the fiber-regularity $\Freg$ as the maximum
over the fibers of their regularity.

\bd In the setting of Definition \ref{RegSch}, 
$$
\tilde{a}_i (\sF ):=\sup_{y\in Y}a^i (\sF \otimes_{\OO_Y}  k(y)), \quad \Freg (\sF ):=\max_{i}\{ \tilde{a}_i (\sF )+i\} ,
$$
and $\Freg (Z):=\max_{i\geq 1}\{ \tilde{a}_i (\sym_Y (\E )/\sI_Z)+i\}$.
\ed

Notice that $\Freg (\sF )$ is finite if $Y$ is covered by finitely many affine charts and $\sF$ is coherent. This holds since the regularity of a graded module over
a polynomial ring over a field is bounded in terms of the number of generators and the degrees of generators and relations
(see {\it e.g.} \cite[3.5]{CFN}).

We now return to the problem of studying the ending degree of local cohomologies of
powers of a graded ideal $I$ in a positively graded Noetherian algebra $A$.

From the comparison of cohomology of stalks and cohomology of fibers, we get from Theorem \ref{regpowgeo} the
following result. As in \ref{regpowgeo} we take a geometric language and do not introduce a graded module (or a sheaf) 
to make the exposition more simple. In case a more general statement is needed, it can be easily derived by
using Theorem \ref{AB2} in place of Theorem \ref{regpowgeo}. The six statements are not independent, but we consider that each of them gives answer to a
question that is quite natural to ask. Notice that  (iv) is essentially equivalent to one of the main results of Eisenbud and Harris  \cite[2.2]{EH}.

\brm
It follows from Theorem \ref{regpowgeo} that the dimension of any fiber of the projection  $\pi$ of the graph to its image (see \ref{regpowgeo}
or below for the precise definition of $\pi$) is bounded above by the cohomological dimension of $A/I$ with respect to $A_+$.
\erm

\bt\label{RegPow}
Let $A:=A_0[x_0,\ldots ,x_n]$ be a positively graded Noetherian algebra and $I$ be a graded $A$-ideal generated
by $m+1$ forms of degree $d$. Set $Y:=\spec (A_0)$ and $X:=\proj (A/I)\subset \proj (A)\subseteq  \tilde{\bf P}_Y^n$. 
Let $\phi :\tilde{\bf P}_Y^n\setminus X\lra {\bf P}^{m}_Y$ be the
corresponding rational map, $W$ be the closure of the image of $\phi$, and
$$
\Gamma \subset \tilde{\bf P}^n_{W}\subseteq \tilde{\bf P}^n_{{\bf P}^{m}_Y}=\tilde{\bf P}_Y^n \times_Y {\bf P}^{m}_Y
$$ be the closure of the graph of $\phi$. Let  
$
\pi : \Gamma \lra W
$ 
be the projection induced by the natural map $ \tilde{\bf P}^n_{{\bf P}^{m}_Y}\lra {\bf P}^{m}_Y$. Then,\\

(i) 
$
\lim_{t\ra +\infty}(\reg ((I^t)^{sat})-dt)
=\max_{i\geq 2}\{ a^i (\Gamma )+i\}.
$
\\

(ii) If $\pi$ admits a fiber $Z\subseteq \tilde{\bf P}_{\spec \ka }^n$ of dimension $i-1$, then 
$$
\lim_{t\ra \infty}(a^{i} (I^t)+i-td)\geq a^{i}(Z)+i=\tilde a^{i}(Z)+i\geq 0.
$$

(iii) Let $\delta$ be the maximal dimension of a fiber of $\pi$. Then,
$$
a^{\delta +1}(I^t) -td=a^{\delta +1}(\Gamma )=\tilde a^{\delta +1}(\Gamma ),\ \forall t\gg 0.
$$

(iv) If $\pi$ is finite, for instance if $X=\emptyset$, then 
$$
\reg (I^t)=a^{1}(I^t)+1=\Freg (\Gamma )+td,\quad \forall t\gg 0
$$
and $\lim_{t\ra \infty}(a^{i}(I^t)-td)=-\infty$ for $i\not= 1$.\medskip

(v)  If $\pi$ has fibers of dimension at most one, for
instance if the canonical map $X\ra Y$ is finite, then
$$
\reg (I^t)-td=\reg (\Gamma )\geq \Freg (\Gamma ),\ \forall t\gg 0,
$$
and $\lim_{t\ra \infty}(a^{i}(I^t)-td)=-\infty$ for $i\geq 2$. 

If furthermore $A$ is standard graded and reduced, $\pi$ has fibers of dimension one, 
all of same degree, then $\Freg (\Gamma )=\reg (\Gamma )$,
$$
\lim_{t\ra \infty}(a^{1}(I^t)-td)\geq \tilde a^{1}(\Gamma)
$$ 
and  equality holds if  $\reg (I^t)=a^{1}(I^t)+1$ for $t\gg 0$.\medskip

(vi)  If $A$ is reduced and, for every connected component $T$ of   $W$, the Hilbert quasi-polynomials of fibers of $\pi$ over any two points 
in $\spec (T)$ differ by a periodic function, then
$$
\reg (I^t)=\Freg (\Gamma )+td,\ \forall \mu\gg 0.
$$

\et

{\it Proof.} (i) is a direct corollary of Theorem \ref{regpowgeo}. Statements (ii), (iii) and (iv) follow from Theorem \ref{regpowgeo} 
and Proposition \ref{CohFib} (i). 

Statements (v) and (vi) follow from  Theorem \ref{regpowgeo}, Proposition \ref{CohFib} (ii) -- notice that $\depth_{S_+}(\R_I )\geq 1$ -- and the equivalence $(i)\Leftrightarrow (iii)$ in 
Proposition \ref{FlatnessCrit} applied on the affine charts covering $\pi (\Gamma )$.\fini


\brm\label{addCEL}
Cutkosky, Ein and Lazarsfeld proved in \cite{CEL} that the limit $s(I):=\lim_{t\ra\infty} \reg ((I^t)^{sat})/t$ exists and is equal to the inverse of a Seshadri constant,  when $A_0$ is a field and $A$ is standard graded. 

Using the existence of $c$ such that $\reg (MI^t)\leq dt+c$ for all $t$ when $I$ is generated in degree at most $d$ and $M$ is finitely generated, one can easily derive the existence of this limit in our more general setting. Indeed, let $r_p:=\reg ((I^p)^{sat})$ and  $d_p:=\min\{ \mu\ \vert\ (I^p)^{sat}=((I^p)^{sat}_{\leq \mu})^{sat}\}$. One has $d_{p+q}\leq d_p +d_q$, hence $s:=\lim_{p\ra\infty} (d_p /p)$ exists.
For any $p$ there exists  $c_p$ such that $\reg (((I^p)^{sat}_{\leq d_p})^tI^q)\leq td_p+c_p$ for all $t\geq 1$ and $0\leq q<p$. The inequalities  $d_{pt+q}\leq r_{pt+q}\leq td_p+c_p$ show that $\lim_{p\ra\infty} (r_p /p)=s$ and that $d_p\geq ps$ for all $p$. 

The same argument applies to any graded ideal $J$ such that $\proj (A/J)\ra Y$ is finite ({\i.e.} $\cd_{A_+}(A/J)\leq 1$). Setting  $r^J_p:=\reg (I^p:_A J^\infty )\leq \reg (I^p)$ and defining $d^J_p$ similarly as above,
$$
d^J_p:=\min\{ \mu\ \vert\ ((I^p:J^\infty)_{\leq\mu}):J^\infty =I^p:J^\infty\},
$$
 the 
limits of $r^J_p/p$ and $d^J_p/p$ exist and are equal. For example, if $X$ is a scheme with isolated 
non locally complete intersection points, then $\lim_{p\ra\infty} \reg (I^{(p)}/p)$ exists, where $I^{(p)}$ denotes the $p$-th symbolic power of $I$.

On the other hand, when $A/J$ has cohomological dimension $2$ it may be that $\reg (I:J^\infty )>\reg (I)$
for $J$ an embedded prime of $I$. This shows that the above argument is not directly applicable for symbolic powers in general. It however implies
that for any $J$, $s^J:=\lim_{p\ra\infty} (d^J_p /p)$ exists and is equal to $\lim_{p\ra\infty} (\rho^J_p /p)$ where 
$$
\rho^J_p:=\min\{ \reg (K)\ \vert\ K\subseteq (I^p:J^\infty ), K:J^\infty =I^p:J^\infty \} .
$$
\erm

\brm
If $I$ is generated in degree at most $d$, Theorem \ref{RegPow}  implies that $s(I)<d$ if and only 
if the morphism $\pi$ corresponding to the ideal $(I_d)$ is finite. More precisely, by Remark \ref{addCEL}, $\pi $ is finite if and only if 
 $\proj (A/I^t)$ is defined by equations of degree $<dt$ for some $t$, and if not $\reg ((I^t)^{sat})-td$
is a non-negative constant for $t\gg0$.

This has been remarked by Niu in \cite{Ni}, using the
definition of $s(I)$ as (the inverse of) a Seshadri constant.
\erm
Theorem \ref{RegPow} also has a consequence on the dimension of the fibers. Assume for simplicity that $A_0$ is a field. 
Set $X:=\proj (A/I)$, with $I$  generated in  degree at most $d$ and let $0\leq i\leq \dim X$.  

Part (ii) in \ref{RegPow} then shows that the morphism $\pi$ associated to  $(I_d )$  has no fiber of dimension greater than $i$ if
there exists $p\geq 1$ and an ideal $K$, generated in degree less than $pd$, such that $\proj (A/I^t)$ and $\proj (A/K)$ coincide locally at each point $x\in {\bf P}^n$ of dimension at least $i$. Indeed if this happens, then
$$
H^j_{A_+}(A/I^{ps})\simeq H^j_{A_+}(A/K^s),\quad \forall j>i,\ \forall s\geq 1,
$$ 
and therefore there exists $c_p$ such that for all $s$ and $j\geq i$, $a^j(I^{ps}) \leq (pd-1)s+c_p$,
showing that $\lim_{t\ra \infty} (a^j(I^{t})-td)=-\infty$ for $j\geq i$.

\end{document}